\documentclass[a4paper]{article}
\usepackage[a4paper]{geometry}
\usepackage{amsthm}
\usepackage{amssymb}
\usepackage{cite}
\usepackage[sc]{mathpazo}
\usepackage{amsmath,scalefnt}
\usepackage{graphicx}
\usepackage{subcaption}
\usepackage{tikz}
\usepackage{mathtools}
\usepackage[authoryear,round]{natbib}
\usepackage{verbatim}
\usepackage{fullpage}
\usepackage[hidelinks]{hyperref}
\allowdisplaybreaks
\usepackage{booktabs}
\usepackage{multirow}
\usepackage{soul}
\usepackage{tikz,pgfplots}

\DeclarePairedDelimiter{\nint}\lfloor\rceil

\newtheorem{example}{Example}

\begin{document}

\title{Using 3D-printing in disaster response: \\ The two-stage stochastic 3D-printing knapsack problem
 }


\author{Denise D. T\"onissen\\
  \footnotesize{
  School of Business and Economics, Vrije Universiteit Amsterdam, Amsterdam, The Netherlands} \\
\footnotesize{d.d.tonissen@vu.nl}
\and
Loe Schlicher\\
\footnotesize{School of Industrial Engineering, Eindhoven University of Technology, Eindhoven, The Netherlands} \\
\footnotesize{l.p.j.schlicher@tue.nl}}
\date{\today}
\maketitle

\noindent \begin{abstract} In this paper, we will shed light on when to pack and use 3D-printers in disaster response operations. For that, we introduce a new type of problem, which we call the two-stage stochastic 3D-printing knapsack problem. We provide a two-stage stochastic programming formulation for this problem, for which both the first and the second stage are NP-hard integer linear programs. We reformulate this formulation to an equivalent integer linear program, which can be efficiently solved by standard solvers. Our numerical results illustrate that for most situations using a 3D-printer is beneficial. Only in extreme circumstances, where the quality of printed items is extremely low, the size of the 3D-printer is extremely large compared to the knapsack size, when there is no time to print the items, or when demand for items is low, packing no 3D-printers is the best option.
\end{abstract}

\bigskip \bigskip

\section{Introduction}
\label{sec:intro}

\noindent Natural disasters, like floods, hurricanes, tornadoes, earthquakes or volcanic eruptions, do have a significant impact on humanity. They cause tens of thousands of deaths, hundreds
of thousands of injuries, and billions of dollars in economic losses each year around the world \citep{dilley2005natural}. The \cite{test2} show that a significant part hereof is caused in the aftermath of such natural disasters. Considering this fact, it is important that the allocation of basic survival resources, such as clean water, food, medicines, cooking utensils, and shelter, as well as the allocation of supportive resources, such as spare parts and repair tools is executed in a fast and proper way \citep{perry2007natural,ozdamar2004emergency}.   \bigskip

\noindent An aspect that complicates the allocation of basic survival and supportive resources is the level of uncertainty. Typically, the impact of a natural disaster is hard to predict \citep{berkes2007understanding,van2006humanitarian} and, as a consequence, the actual number of basic survival and supportive resources needed, is hard to predict as well. A new, emerging technique that can deal with this form of uncertainty is the technique of 3D-printing, (see, e.g., \citet{savonen2018development,tatham2015three} and the references therein).
With this technique, also known as additive manufacturing, one can create objects (e.g., cooking utensils, bottles, screwdrivers, shelter, pipe clamps, or water purification kits) locally. This means that, if one decides to transport 3D-printers to a disaster area, one introduces the possibility to anticipate on actual (i.e., local) demand. For instance, in Nepal, this technique has been successfully used for printing medical supplies and restoring water supply in the aftermath of the earthquake in 2015 \citep{maldives,maldives3} and in Ha\"{i}ti, this technique has been successfully used for printing screwdrivers, pipe clamps and bottles in the aftermath of the earthquake in 2010 \citep{maldives2}. \bigskip

\noindent The benefits of using 3D-printing are not limited to the possibility to deal with demand uncertainty.
For instance, in literature, it is believed that printing material, which is in liquid or powder form,
can be packed much more efficient than physical items.
Consequently, more printing material (and subsequently more (future) printed items) can be packed for a disaster response mission. However, 3D-printing also has some drawbacks, which makes it a less attractive option. For instance, the weight and volume of 3D-printers is still significant nowadays. This implies that, if one plans to transport some 3D-printers to a disaster area, this weight and volume cannot be used to pack other critical items. Moreover, the quality, and usability, of printed items is still, inferior to non-printed items. And, finally, \textcolor{black}{the 3D-printing time of items} is still in terms of hours --time that is crucial during a disaster response mission. Indeed, with these advantages and disadvantages of 3D-printing (see Table 1 for a summary), is not clear immediately whether we should pack 3D-printers for a disaster response mission, or not.\bigskip

\begin{table}[h!] \centering
\begin{tabular}{cc} \hline
Positive effects & Negative effects \\ \hline
\emph{Possibility to deal with demand uncertainty} & \textcolor{black}{\emph{3D-printing time}} \\ 
\emph{Efficient packing of printing material} & \emph{Weight and volume of 3D-printer\textcolor{black}{s}}\\
 & \emph{Reduced quality \textcolor{black}{of} printed item\textcolor{black}{s}} \\   \hline
\end{tabular} \caption{Positive and negative effects of using 3D-printing.}
\label{tab:3Dposneg}
\end{table}

\noindent In this paper, we will shed light on when to bring 3D-printers or not, which amount of printing material to take and which items to bring physically, to a disaster area. For that, we introduce a new type of problem, which we call the two-stage stochastic 3D-Printing knapsack problem (TSS-3DKP). In the first stage of this problem, a decision maker has to fill a multidimensional knapsack (e.g., a cargo airplane at a home location), with the special consideration of taking 3D-printers with an associated amount of printing material, or not.
Then, in the second stage (e.g., in the disaster area), demand for basic survival requirements as well as supportive resources is revealed, and the physically brought items that match demand are allocated.
If at least one 3D-printer is taken, the next consideration is how to use the printing material for the remaining demand.
Based on a maximum of printing time per 3D-printer  - which resembles the urgency of delivering items fast (e.g., within 24 or 48 hours) - a decision on the number and type of printed items (per 3D-printer) is made.  Meeting demand, via a physically brought item or a printed item, results in an item-specific reward, which may depend on its nature (i.e., whether it is printed or not). This reward may represent the (relative) importance of having such an item. Aim is to make, a priori, a decision upon the number of physical items to take, whether to bring 3D-printers or not, the amount of printing material to take, and how to use it (i.e., which items to print), in order to maximize the expected total reward. \bigskip

\noindent From a practical perspective, it is important that our
TSS-3DKP can be solved within reasonable computation time (e.g., within hours). This may be challenging, since our problem is, in contrast to
standard two-stage stochastic programming problems \citep{Birge:2011}, NP-hard in both the first and the second-stage
problem.
We tackle our TSS-3DKP, by transforming it into a deterministic equivalent (i.e., into a large integer linear programming problem).
This allows us to solve our problem by standard solvers within reasonable computation time.
We want to emphasize that this transformation is not straightforward,
because the number of constraints and variables of our second-stage problem depend on the first-stage decisions.
We overcome this dependency by introducing several dummy variables  and a smart upper bound on the possible number of packed 3D-printers.
\bigskip

\noindent We use \textcolor{black}{our} deterministic equivalent to execute several numerical experiments. In particular, we will investigate how the advantages and disadvantages of Table 1 (i.e., quality of printed items, the weight and volume of 3D-printers, the storage efficiency of printing material, the printing time, and the demand uncertainty) effect the decision of when to use 3D-printing, and when not. \bigskip

\noindent Now, we summarize the main contributions of this paper:

\begin{enumerate}
\item{We formulate a new two-stage stochastic knapsack problem, inspired by disaster
response  operations, in which one determines the numbers of physical items and 3D-printers
to pack, the amount of printing material to pack, and how to use this printing
material.}
\item{We are able to provide a deterministic equivalent for our two-stage stochastic 3D-printing knapsack problem, which allows us to solve the problem within reasonable computation time.}
\item{Via numerical experiments, we are able to identify under which circumstances it is beneficial to bring 3D-printers, and by how much, to a disaster area.}
\end{enumerate}



\noindent Finally, it is worth mentioning that our model formulation of the TSS-3DKP and associated (numerical) results may find future applications in other domains such as the military (e.g, a military mission on a remote location) and aerospace (e.g., a long-term mission to Mars). \bigskip

\noindent The remainder of this paper is organized as follows. We start in Section 2 with a literature review. In Section 3, we introduce our two-stage stochastic programming problem. In Section 4, we focus on solving our problem, by presenting and testing a deterministic equivalent. Then, in Section 5, we identify, via numerical experiments, under which circumstances it is beneficial to bring 3D-printers, and by how much, to a disaster area. We close with a conclusion in Section 6.

\section{Literature review}

\noindent In this section, we will identify how our TSS-3DKP contributes to the literature. We do so by providing an overview of the (OM/OR) literature on two-stage stochastic knapsack problems, disaster operations management, and 3D-printing, and for each of them, identify how our TSS-3DKP contributes to it.

\subsection{Two-stage stochastic knapsack problems}

\noindent A classical combinatorial optimization problem is the knapsack problem (\citet{Kolesar:1967}).
In this problem, a decision maker has to fill a knapsack with items, each with an associated weight and reward.
The aim of the decision maker is to find the best layout, i.e., a feasible combination of items, according to the size of the knapsack, that maximizes the sum of the rewards.
Many variations of this classical problem have been studied.
For some of these variations, the decision maker has to fill a knapsack based on incomplete information
(e.g., the reward of items is uncertain) and has recovery options once complete information becomes available.
The aim of the decision maker is then to maximize the expected reward.
In literature, these variations are classified as two-stage stochastic knapsack problems.
The first stage refers to the setting with partial information and the second stage refers to the setting with complete information and recovery options. \bigskip

\noindent The literature on two-stage stochastic knapsack problems is relative young and consists of a few papers only.
\cite{Kosuch:2011} study a two-stage stochastic knapsack problem with normally distributed item weights and a recovery option that limits to either the addition or the removal of items.
To restrict the probability that the item weights exceed the knapsack size (in the second stage) they include a chance constraint to the first-stage problem.
A variant, with discretely distributed item weights, is studied by \cite{Kosuch:2014}.
Knapsack size uncertainty is studied by \cite{Akker:2016}.
They model this uncertainty by discrete scenarios and introduce a recovery option that consists of the removal of items.
This problem is extended to multiple knapsacks by \citet{Tonissen:2017}.
Finally, a quadratic version of the two-stage stochastic knapsack problem with item weight uncertainty has been studied by \citet{Lisser:2010}
and with both item weight and reward uncertainty by \citet{Song:2018}.\bigskip

\noindent Our paper contributes to this new stream of literature, by studying a two-stage stochastic knapsack problem with uncertainty in (item) demand and a recovery option that consists of 'creating' new items out of another item (i.e., the printing material). To the best of our knowledge, this paper is the first that studies a two-stage stochastic knapsack problem with such type of uncertainty and recovery option.

\subsection{Disaster operations management}

Disaster operations represent the set of activities performed before, during and after a disaster in order to diminish its impact (\citet{altay2006or}).
In general, it is hard to prepare for, and consequently to manage such type of activities. This is mainly due to the unpredictable nature of disasters. In the OM/OR literature on disaster operations management, this uncertainty plays a predominant role (see, e.g., the extensive literature reviews of \citet{Galindo:2013} and \citet{Leiras:2014}). Common uncertainties in this literature include demand uncertainty \citep{Verma:2015,Tricoire:2012,Alem:2016,Tofighi:2016},
uncertainty of the condition of facilities and/or roads \citep{Hu:2019,Fan:2010,Sanci:2019,Li:2012,Chang:2007}, travel time uncertainty \citep{Ahmadi:2015,Doyen:2012,Bayram:2018}, and cost uncertainty \citep{Bozorgi:2012,Li:2011,Paul:2019}. \bigskip

\noindent  Two-stage stochastic programming is an appropriate method to deal with such uncertainties.
Instead of using standard cost minimizing objectives, in disaster operations management, meeting expected demand is often preferred due to significant consequences with respect to human lives (\citet{Barbarosoglu:2004} and \citet{Rawls:2010}). In \citet{Balcik:2008} the fulfillment of demand is modelled by maximizing the expected satisfied demand, in   \citet{Salmeron:2010} by minimizing the expected casualties, and in \citet{Noyan:2016} by maximizing the expected accessibility. A common first-stage decisions is the in-advance storage of relief items (\citet{Davis:2013,Lodree:2012}) or locating facilities (\citet{Elcci:2018,Li:2011}). A few authors use other first-stage decisions such as the retrofitting of roads (\citet{Peeta:2010}), buildings (\citet{Zolfaghari:2015}) or bridges (\citet{Liu:2009}). The most prevalent second-stage decisions are the transport of commodities in the aftermath of a disaster (\citet{Rezaei:2016,Tofighi:2016}) or an evacuation plan (\citet{Li:2011,Li:2012}). We refer to \citet{Grass:2016} for an overview of two-stage stochastic programming problems for disaster operations management. \bigskip

\noindent In this paper, we also use two-stage stochastic programming to deal with uncertainty (namely, demand uncertainty). In particular, we will use this method in a disaster response setting.
To the best of our knowledge, we are the first who use two-stage stochastic programming to investigate whether bringing 3D-printers to a disaster area is useful or not.
We want to emphasize that there exist some qualitative papers on 3D-printing for disaster response missions (see, e.g., \citet{savonen2018development} and \citet{Rodriguez:2018}). These studies focus on the necessary requirements and specifications 3D-printers should have for being successfully used in a disaster area.

\subsection{3D-printing}

\noindent Although the technique of 3D-printing, also known as additive manufacturing, has already been applied in practice for years, the development of quantitative models studying the impact of 3D-printing in OM/OR literature is rather limited. To the best of our knowledge, there exist a few published works in this domain only. \citet{westerweel2018traditional} investigate the impact of 3D-printing on component design. In particular, they characterise under which conditions a component should be produced with traditional technology and under which conditions a component should be produced via 3D-printers. \citet{song2019stock} present a general framework to study the design of spare parts logistics in the presence of 3D-printing technology. In particular, they formulate a model that determines which parts to stock and which to print. \citet{khajavi2014additive} and \citet{liu2014impact}
 both model a multi-echelon spare parts supply chain and numerically investigate the
effect of centralized versus decentralized 3D-printing capacity. \citet{dong20163d} evaluate the choice between traditional technology and 3D-printing related to assortment planning for general inventories in a manufacturing setting. An overlap in all these papers is the presence of a 3D-printer with a, possibly unlimited, amount of printing material. However, in our problem, the presence of a 3D-printer and associated printing material is not guaranteed: it is an essential decision in our optimization problem. This indicates that we do study a new type of 3D-printing (OM/OR) problem and so contribute to this rather new stream of literature.

\section{Two-stage stochastic 3D-printer knapsack problem}
\label{sec:tsp}

\noindent In this section, the two-stage stochastic 3D-printing knapsack problem (shortly TSS-3DKP) will be formulated. We will do so by first describing the two stages (i.e., the first and second stage) and thereafter presenting the mathematical formulation of our TSS-3DKP.

\subsection{The first stage}

\noindent In the first stage, a decision maker has to fill a multidimensional knapsack (e.g., a cargo airplane at a home location) with physical items, units of printing material and 3D-printers such that it does not exceed the weight capacity $W$ and volume capacity $V$ of the knapsack. Each physical item $i$, from the set of physical items $N$, has a weight $w_{i}$ and a volume $v_{i}$. In addition, each unit of printing material has weight $w_b$ and volume $v_b$ and each 3D-printer has weight $w_p$ and volume $v_p$. The \textbf{first-stage decision of the decision maker} is denoted by vector $\textbf{a} = ((a_i)_{i \in N},a_p,a_b)$ with $a_i$ the number of times item $i \in N$ is added to the knapsack, $a_b$ the units of printing material included into the knapsack, and $a_p$ the number of 3D-printers added to the knapsack.

\subsection{The second stage}

\noindent In the second stage (e.g., a disaster area) the decision maker has to allocate the $a_b$ units of printing material to the $a_p$ 3D-printers.
Before the decision maker does so, demand for items is revealed and the physically brought items that match (this) demand are allocated.
Demand is modelled by a discrete set of scenarios $S$,
where each scenario $s \in S$ specifies the demand $d_i^s$ for all items $i \in N$ and occurs with probability $q_s \in [0,1]$ such that $\sum_{s \in S} q_s = 1$. We want to stress that, in literature on disaster operations management, discrete scenarios are often used to model uncertainty. This is due to the uniqueness of disasters, which makes the determination of probability distributions problematic (see, e.g., \citet{Grass:2016}). For each scenario $s \in S$, we denote the number of physically brought items that match demand by $a_i^s = \min\{a_i,d_i^s\}$ for all $i \in N$ and
 consequently, denote remaining demand by $d_i^s-a_i^s$ for all $i \in N$.
Based on this remaining demand, the decision maker has to allocate the $a_b$ units of printing material to the $a_p$ 3D-printers.
For this second-stage decision, the decision maker has to take into account the maximum amount of printing time $T$ per 3D-printer,
and the fact that each printable item $i \in N^p$, with $N^{p} \subseteq N$, requires $m_i$ units of printing material and time $t_i$ to print.
Finally, given that $P(x) = \{i \in \mathbb{N}_{+} \vert i \leq x \}$ is defined as the set of 3D-printers for any $x\in \textcolor{black}{\mathbb{N}_{0}}$, the \textbf{second-stage decision of the decision maker} for each scenario $s \in S$ is denoted by $\textbf{p}^s=(p_{ij}^s)_{i \in N, j \in P(a_p)}$, with $p_{ij}^s$ the number of times item $i$ is printed on 3D-printer $j$ in scenario $s$. \bigskip

\noindent Meeting demand via a physically brought item $i \in N$ or printed item $i \in N^p$ results in an item-specific reward,
which may depend on its nature (i.e., whether it is printed or not).
This reward may represent the (relative) importance of having such an item.
The reward for a physically brought and matched item $i \in N$ is $r_i$ and its printed reward is $\alpha  r_i$ with $\alpha \in [0,1]$ for all $i \in N^p$.
 This factor $\alpha$ resembles the lower quality of the printed items.
The total reward for a given first-stage decision $\textbf{a}$ and a second-stage decision $\textbf{p}^s$
for scenario $s\in S$ is given by $\sum_{i \in N} a_i^s r_i + \alpha \sum_{i \in N^p} \sum_{j \in P(a_p)}p_{ij}^s r_i$.
And, subsequently, the expected total reward, based on all possible demand scenarios, is given by $\sum_{s \in S} q_s \left[ \sum_{i \in N} a_i^s r_i + \alpha \sum_{ \textcolor{black}{i \in N^p}} \sum_{j \in P(a_p)}p_{ij}^s r_i\right]$.
Aim of the decision maker is to make a first-stage decision $\textbf{a}$ and a second-stage decision  $\textbf{p}^s$ for all $s \in S$ that maximizes the expected total reward\footnote{For the specific setting with $r_i=1$ for all $i \in N$, our objective can be recognized as the (expected) amount of satisfied demand. This objective is used frequently in the literature on disaster operations management.}.

\subsection{Mathematical formulation of \textcolor{black}{the} TSS-3DKP}

\noindent Now we give a mathematical formulation of our two-stage stochastic 3D-printing knapsack problem.

\begin{displaymath} \mbox{TSS-3DKP} \quad \max  \sum_{s \in S} q_s \text{Q}(\boldsymbol{a},s)\end{displaymath}
 \begin{align}\text{s.t. }
&& \sum_{i \in N} a_{i} w_{i} + a_p w_p + a_b w_b  \le W,  \label{1.1} \\
&& \sum_{i \in N} a_{i} v_{i} + a_p v_p + a_b v_b  \le V, \label{1.2} \\
&&  a_b \le a_p \mathcal{M} \label{1.3},  \\
&& a_p, a_b, a_{i} \in \mathbb{N}_{0} && \forall i \in N, \label{1.4}
 \end{align}

\noindent where,

\begin{displaymath} \text{Q}(\boldsymbol{a},s):=  \quad  \max \sum_{i\in N} a_{i}^{s} r_{i} + \sum_{i \in N^p}  \sum_{j\in P(a_p)} \alpha p_{ij}^{s}  r_{i}  \end{displaymath}
 \begin{align}\text{s.t. }
&& a_i^{s} \le d^{s}_{i} && \forall i \in N \setminus N^p, \label{2.0}\\
&& a_i^{s} + \sum_{j\in P(a_p)} p_{ij}^{s} \le d^{s}_{i} && \forall i \in N^p, \label{2.1}\\
&& a_i^{s} \le a_{i} && \forall i \in N, \label{2.2}\\
&& \sum_{i \in N^p} \sum_{j \in P(a_p)} p_{ij}^{s} m_{i}  \le a_b,  \label{2.3}  \\
&& \sum_{i \in N^p} p_{ij}^{s} t_{i} \le T && \forall j \in P(a_p), \label{2.4}\\
&& a_{i}^{s} \in  \mathbb{N}_{0} &&  \forall i \in N, \label{2.5}\\
&& p_{ij}^{s} \in  \mathbb{N}_{0} && \forall i \in N^p, \quad \forall j \in P(a_p) \label{2.6}.
 \end{align}

\noindent Constraints \eqref{1.1} and \eqref{1.2} ensure that no more number of items, units of printing material and 3D-printers will be packed, according to the \textcolor{black}{weight and volume} restrictions of the knapsack. Constraint \eqref{1.3}, which is a big-$\mathcal{M}$ formulation, guarantees that it is only allowed to pack printing material if at least one 3D-printer is taken. Constraint \eqref{1.4} ensures integrality of the decision variables. Constraints \eqref{2.0}, \eqref{2.1}, and \eqref{2.2} ensure that reward can only be gained for those physical items and printed items that meet demand. Constraint \eqref{2.3} guarantees that the amount of printing material used does not exceed the amount of printing material packed.
Constraints \eqref{2.4} ensure that each 3D-printer can be used for at most $T$ time units.
Constraints \eqref{2.5} and \eqref{2.6} ensure integrality of the decision variables. Note that we model $a_i^s(=\min\{a_i,d_i^s\})$ as a decision variable in our second-stage problem which takes, due to our formulation, the minimum of the values $a_i$ and $d_i^s$ for all $i \in N$ and all $s \in S$.
Also note that the number of constraints \eqref{2.4} and \eqref{2.6} depend on the first-stage decision $a_p$. Hence, the number of constraints and variables of the second-stage problem depend on the decisions in the first-stage problem.\bigskip

\noindent  We want to emphasize that both the first-stage and the second-stage problem of our TSS-3DKP are NP-hard integer programming problems. Both problems are NP-hard because they are generalisations of the unbounded knapsack problem \citep{Martello:1990}, which is a classical NP-hard combinatorial optimization problem. We will support these claims, by showing that special cases of our first and second-stage problems can be recognized as unbounded knapsack problems. For the first-stage, we consider the special case with no 3D-printers, no printing material, and exactly one scenario for which demand for each item is infinite. Furthermore, we set $V=0$ and $v_{i} = 0$ for all $i \in N$.
For this situation, our first-stage problem can be recognized as a classical unbounded knapsack with size (i.e., weight capacity) $W$ and items with weight $w_i$ and reward $r_{i}$ for all $i \in N$.
For the second-stage, we consider the special case with no physical items, exactly one 3D-printer, infinite demand, and $m_i=0$ for all $i \in N^p$.
For this situation, our second-stage problem can be recognized as a classical unbounded knapsack with size $T$ and items with weight $t_i$ and reward $\alpha r_{i}$ for all $i \in N^{p}$. \bigskip

 \noindent  We will now illustrate our two-stage stochastic 3D-printing problem by means of a small example.

\begin{example} Consider a knapsack with \textcolor{black}{weight capacity $W=4$ and volume capacity $V=4$,}
that can be filled with two types of printable items (i.e., $N=N^{\textcolor{black}{p}}=\{1,2\}$).
For the items, we have \textcolor{black}{$w_1=v_1=w_2=v_2=4$} with $r_1=1$ and $r_2=2$.
\textcolor{black}{The items require $m_1 = m_2 = 2$ units of printing material} and each unit of printing material has weight $w_b=1$ and volume $v_b=1$.
In addition, each 3D-printer has weight \textcolor{black}{$w_p=2$} and volume \textcolor{black}{$v_p=2$}, and,
in total, there are two scenarios (i.e., $S= \{s_1,s_2\}$), with $q_{s_1}=0.7$ and $q_{s_2}=0.3$.
Demand is given by $d^{s_1}_1 = 1, d^{s_1}_2=0$, $d^{s_2}_1= 0, d^{s_2}_2= 1$.
Finally, the maximum amount of printing time per 3D-printer is given by $T=1$, the printing time per item by $t_1=t_2=1$, and the quality factor (per printed item) by $\alpha =0.8$.\bigskip

\noindent For this specific situation, there exist three candidates solutions, namely: $(i)$
taking one physical item of type 1, ($ii$) taking one physical item of type 2, or ($iii$) taking one 3D-printer with \textcolor{black}{2 units} of printing material.
 For the first two candidate solutions, no 3D-printer is taken, and consequently, the expected total reward for the first strategy equals $(0.7 \cdot 1=) 0.7$
 and for the second strategy equals $(0.3 \cdot 2=) 0.6$.
 For determining the expected reward for the third strategy, we first need to identify how the units of printing material should be used in the two demand scenarios.
  For the first demand scenario (i.e., $s_1$), it is optimal to print item 1, resulting in $(0.7 \cdot 0.8 \cdot 1=)0.56$,
  and for the second demand scenario (i.e., $s_2$), it is optimal to print item 2, resulting in $(0.3 \cdot 0.8 \cdot 2 =)0.48$.
  Hence, the expected total reward under the third strategy equals $(0.56+0.48=)1.04$, implying that bringing a 3D-printer and allocating the printing material according to demand is optimal.
 \end{example}

\noindent To the best of our knowledge, there exist no standard (commercial) solvers that can directly solve two-stage stochastic programming problems.
For that reason, we will transform our TSS-3DKP into another, but equivalent, formulation that is suitable for standard solvers (e.g., CPLEX).

 \section{Solving the TSS-3DKP}
\label{sec:ILP}

In this section we will present an equivalent formulation of our TSS-3DKP that can be solved by standard solvers. Moreover we will show, via some numerical experiments,
that this equivalent formulation can be solved (with standard solvers) within reasonable computation time.

\subsection{An equivalent formulation for solving the TSS-3DKP}

It is well-known that a standard two-stage stochastic programming problem can be modelled as a large integer linear programming problem with
variables and constraints for each of the scenarios.
Such an integer linear programming formulation is called the deterministic equivalent. The idea is to transform our TSS-3DKP, which is also a two-stage stochastic programming problem, into its deterministic equivalent. However, the TSS-3DKP has the non-standard feature that the number of constraints \textcolor{black}{and variables} in the second-stage depends on the first-stage decisions, implying that the deterministic equivalent cannot be given directly.
We overcome this dependency, by identifying an upper bound on the total number of 3D-printers that can be \textcolor{black}{packed}, and include printed-related constraints \textcolor{black}{and variables} of the second-stage problem into our deterministic equivalent as if the knapsack would be filled with this upper bound of 3D-printers. Although this approach may lead to an increase in the total number of constraints and variables, it allows for a  reformulation of our TSS-3DKP\footnote{In Section 4.2, we will also show that, although the total number of constraints \textcolor{black}{and variables} grows under this approach, our deterministic equivalent formulation can still be solved in reasonable computation time.}. \bigskip

 \noindent Before we present the deterministic equivalent formulation, we first explain how we determine an upper bound on the maximum number of \textcolor{black}{packed} 3D-printers.
For each scenario $s \in S$, we determine how many 3D-printers are needed to print all demand $d^s$, \textcolor{black}{assuming that items are allocated one by one in order of their indices to the 3D-printers. Hence, we first try to allocate all demand of item 1 to the first 3D-printer, followed by item $2,3,\ldots,|N^p|$, and we go to the next 3D-printer if adding another item would exceed the printing $T$ of the 3D-printer.} Taking the \textcolor{black}{maximum number of 3D-printers} over all possible scenarios then gives us our first upper bound $U$. Sometimes, this \textcolor{black}{upper bound} exceeds the total knapsack \textcolor{black}{weight} or volume. In those cases, we select the maximum number of 3D-printers that fit the knapsack as upper bound. Formally, we define our \textcolor{black}{upper bound} as $Z = \min{\{\lfloor \frac{W}{w_p} \rfloor, \lfloor \frac{V}{v_p} \rfloor, U\}}$. \bigskip

\noindent We will now illustrate our upper bound by means of a small example.

\begin{example} Consider a situation with two types of printable items (i.e., $N=N_p=\{1,2,\}$), two scenarios (i.e., $S=\{s_1,s_2\}$) with associated demands $d^{s_1}_1=3, d^{s_1}_2=1$, $d^{s_2}_1=1, d^{2_2}_2=2$, a knapsack with weight capacity $W = 10$ and volume capacity $V=12$, a 3D-printer with $w_p=5$, $v_p=3$, $T=5$, and printing times $t_1=2$ and $t_2=3$.
\textcolor{black}{We now allocate the demand of the items to the 3D-printers in order of their indices.}
Hence, for scenario $s_1$, we allocate two items of type 1 to the first 3D-printer, and one item of type 1 and one item of type 2 to the second 3D-printer.
For scenario $s_2$, we allocate one item of type 1 and one item of type 2 to the first 3D-printer, one item of type 2 to the second 3D-printer, and one item of type 2 to the third 3D-printer.
 The maximum number of 3D-printers needed (over all scenarios) is thus $U=3$.
 Moreover, $\lfloor \frac{W}{w_p} \rfloor = 2$ and $\lfloor \frac{V}{v_p} \rfloor = 4$, implying an
  upper bound on the total number of 3D-printers of $Z=2$.  \end{example}

\noindent For the deterministic equivalent formulation, we also need to introduce some new notation. For the given upper bound $Z$, we define binary variables $y_j \in \{0,1\}$ for all $j \in P(Z)$, with $y_j=1$ \textcolor{black}{if 3D-printer $j$ is packed and 0 otherwise.} \bigskip

\noindent Now, we are ready to present the deterministic equivalent of our TSS-3DKP.

\begin{displaymath}  \text{ILP-3DKP}(Z) \quad \max  \sum_{s \in S} q_s  \left[\sum_{i \in N} a_{i}^{s} r_{i} + \sum_{i \in N^p}  \sum_{j \in  P(Z)} \alpha p_{ij}^{s} r_{i} \right] \end{displaymath}
 \begin{align}\text{s.t. }
& \sum_{i \in N} a_{i} w_{i} + a_b w_b + \sum_{j \in  P(Z)} y_j w_p  \le W,   \label{DET1} \\
& \sum_{i \in N} a_{i} v_{i} + a_b v_b + \sum_{j \in  P(Z)} y_j v_p   \le V,  \label{DET2} \\
& a_b \le \sum_{j \in  P(Z)} y_j M, \label{DET3}  \\
& a_i^{s} \le d^{s}_{i} && \forall i \in N \setminus N^p, \quad \forall s \in S, \label{DET4a}\\
& a_i^{s} + \sum_{j \in  P(Z)} p_{ij}^{s} \le d^{s}_{i} && \forall i \in N^p, \quad \forall s \in S, \label{DET4b}\\
& a_i^{s} \le a_{i} && \forall i \in N, \quad \forall s \in S, \label{DET5}\\
& \sum_{i \in N^p} \sum_{j \in  P(Z)} p_{ij}^{s} m_{i}  \le a_b && \forall s \in S, \label{DET6}  \\
& \sum_{i \in N^p} p_{ij}^{s} t_{i} \le T y_{j} && \forall j \in P(Z), \quad \forall s \in S, \label{DET7}\\
& y_{j} \le y_{j-1} && \forall j \in  P(Z) \setminus \{1\},  \label{DET8} \\
&y_{j} \in \{0,1\} && \forall j \in  P(Z), \label{DET9}\\
& a_b \in \mathbb{N}_{0} && \label{DET10}\\
&a_{i} \in  \mathbb{N}_{0} && \forall i \in N, \label{DET11}\\
& a_{i}^{s} \in  \mathbb{N}_{0} &&  \forall i \in N, \label{DET12}\\
& p_{ij}^{s} \in  \mathbb{N}_{0} && \forall i \in N^p, \quad \forall j \in  P(Z), \quad \forall s \in S. \label{DET14}
 \end{align}

\noindent The objective coincides with the original objective of the TSS-3DKP: it is the expected total reward obtained from
\textcolor{black}{meeting demand by physical and 3D-printed items} over all possible scenarios. Constraints \eqref{DET1} and \eqref{DET2} limit the total item, 3D-printer and printing material weight and volume to the weight and volume \textcolor{black}{restrictions} of the knapsack. The big-$\mathcal{M}$ Constraint \eqref{DET3} guarantees that printing material is only taken when at least one 3D-printer is taken.
This constraint guarantees indirectly that we can only print when we have a \textcolor{black}{3D-}printer.
Constraints \eqref{DET4a}, \eqref{DET4b} and \eqref{DET5} ensure that reward can only be gained for those physical items and \textcolor{black}{3D} printed items that meet demand. Constraints \eqref{DET6} guarantee that the amount of printing material used does not exceed the amount of printing material packed.
Constraints \eqref{DET7} ensure that each \textcolor{black}{3D-}printer can be used for at most $T$ time units.
Note that this amount (of time) becomes zero if the 3D-printer is not taken.
Constraints \eqref{DET8} guarantee that \textcolor{black}{3D-}printer $j$ is taken if and only if 3D-printer $j-1$ is taken.
This excludes many symmetric solutions and consequently decreases the solution time of the CPLEX solver.
Constraints \eqref{DET10}, \eqref{DET11}, \eqref{DET12}, and \eqref{DET14} ensure integrality of the decision variables.

\subsection{Computational experiments}

\noindent \textcolor{black}{In this section, we will discuss the outcomes of our \textcolor{black}{computational} experiments to provide insights in the solvability of our deterministic equivalent formulation. We do so by first explaining how we generate our instances and then continue by discussing our \textcolor{black}{computational} experiments.}

\subsubsection{Instance generation}

\noindent For every instance, we first specify the following parameters: the number of items $|N|$, an upper limit to demand $D$, the number of scenario $|S|$, the (printed item) quality factor $\alpha$, and the 3D-printer weight $w_p$ and volume $v_p$. Thereafter, we will generate the items.
Remember that every item has a weight $w_i$, a volume $v_i$, a reward $r_i$, a required printing material $m_i$, and printing time $t_i$.
We generate the weight $w_i$ and the reward $r_i$ according to the uncorrelated instance class of Pisinger \citep{pisinger:2005}.
This class is realistic for disaster response missions: it considers many different items for which the weight, volume, and rewards vary heavily.
For the volume, printing time and required printing material per item, we use the following distributions for generation: $v_i \sim \nint{U[0.2,5] w_i}$, $t_i \sim \nint{U[0,10]}$ and $m_i \sim \nint{U[0.5, 0.9] \cdot \min{(w_i, v_{i})}}$.
Note that we use $\nint{~}$ to denote that we round our values to the nearest integer.
Moreover, we set $w_b=v_b=1$. In combination with the generation of $m_i$, this implies that printed items \textcolor{black}{can always be
packed more efficiently} than physical items. Thereafter, we generate the demand for the $\vert S \vert$ scenarios.
First, we generate the maximum demand for each item $i \in N$ from uniform distribution $U_i \sim \nint{U[1,D]}$.
Then, for each scenario $s\in S$ we generate a demand $d_{i}^{s} \sim \nint{U[0,U_i]}$ for all $i \in N$ and the probability $q_s = \frac{1}{\vert S \vert}$.
We continue by generating printing
time $T \sim \nint{U[0.2, 1] \cdot \sum_{s \in S}q_s \sum_{i \in N^p}  t_{i} d_{i}^s}$.
Note that this printing time $T$ is a fraction of the average total time to print all items.
Finally, we generate the knapsack weight \textcolor{black}{capacity} $W \sim \nint{U[0.5,1] \cdot \sum_{s \in S}q_s \sum_{i \in N}  w_{i} d_{i}^s}$, which can be seen as a fraction of the average total weight of all items, and volume \textcolor{black}{capacity} $V \sim \nint{U[0.5,2] \cdot W}$, which depends on the knapsack weight capacity.\bigskip

\noindent We will study various instance sets with varying number of items, upper limit demand and number of scenarios, but with the same $\alpha=0.8$ and $w_p=v_p=5000$. We will refer to such an instance set by $NxDySz$ with $x$ the number of items, $y$ the upper limit to the demand and $z$ the number of scenarios.

\subsubsection{Computational results}
\label{sec:comp}
\noindent The experiments are programmed in Java with the CPLEX library version 12.8.0, and run on a laptop with an Intel Core i7-4710MQ Quad Core 2.5 GHz processor with 32 GB of RAM. We use CPLEX standard settings and the reported solution times include everything that is required
to solve our instances, such as the time required for reading the instances, the time to determine upper bound $Z$, the time required for building the model,
and the time to solve the model. The maximum computation time to solve a single instance is set equal to 3600 seconds. \bigskip

\noindent We will test our ILP-3DKP(Z) on five different instance sets, each consisting of 100 instances. The first instance set, which we use as a base case, consists of 100 items, an upper limit of demand of 100, and 50 scenarios. Recall that we (can) refer to this instance set as N100D100S50. In addition, we will study instance set N200D100S50, N100D200S50, N100D100S100, and N200D200S100. \bigskip

\noindent We will study each instance set on four performance criteria: the average number of cuts and nodes used by the branch-and-bound tree of CPLEX, the number of instances that cannot be solved within 3600 seconds, and the average computation time of those instances that can be solved within 3600 seconds.  For the first two criteria (i.e., cuts and nodes) we want to stress that if an instance is not solved within 3600 seconds, a fail is registered, and the number of cuts and nodes used so far is reported. The performance criteria of the five instance sets can be found in Table \ref{tab:ILP}.\bigskip

\begin{table}[h!]
\begin{center}
\scalebox{1.0}{
\begin{tabular}{ | c || c | c |c | c| }
\hline
\textbf{Instance set} & \textbf{cuts (\#)} & \textbf{nodes (\#)}  & \textbf{fails (\#)} & \textbf{time (seconds)}  \\
\hline
N100D100S50 & 1297 & 38822	& 1 &  12 \\ \hline
N200D100S50 & 2870 & 0	& 0 & 3  \\ \hline
N100D200S50 & 1488 & 1	& 0 & 1  \\ \hline
N100D100S100 & 2638 & 19921	 & 1 & 50  \\ \hline
N200D200S100 & 5616 & 5	 & 0 & 16  \\ \hline
\end{tabular}
} 	
\end{center}
\caption{The performance criteria of the 5 instance sets.}
\label{tab:ILP}
\end{table}

\noindent Table \ref{tab:ILP} indicates that almost all instances (except for the two failed ones) can be solved within minutes. In particular, we observe, based on individual results not reported here, that most of these instances do not require branching (i.e., generate no nodes) and are solved entirely by adding cuts.
However, the two failed instances (one from N100D100S50 and one from N100D100S100) turn out to be hard to solve (i.e., they cannot be solved to optimality within 3600 seconds). These particular instances have a high number of generated nodes (namely 3017763 and 1158576) and consequently require extensive branching.
Moreover, it turns out that, if we increase the maximum computation time to 86,400 seconds (i.e., 24 hours) these two instances cannot be solved to optimality either. \bigskip

\noindent
Since finding an optimal solution seems to be a problem for some instances, the next step is to investigate the effect of the stopping criteria on the branch-and-cut algorithm (in CPLEX). The standard stopping criteria of this algorithm is the relative gap of 0.01\%. We will study the effect by solving a specific instance set for various relative gaps, namely 0.01\%, 0.1\%, 1.0\%, and 10\%. The specific instance set combines the five, previously introduced, instance sets. We present our results in Table \ref{tab:ILPgap}.  \bigskip

\begin{table}[h!]
\begin{center}
\scalebox{1.0}{

\begin{tabular}{ | l || c |  c | c| c| }
\hline
\textbf{Relative gap} & \textbf{cuts (\#)} & \textbf{nodes (\#)}  & \textbf{fails (\#)} & \textbf{time (seconds)}  \\
\hline
0.01\%  & 2782  & 11750 & 2  & 31 \\ \hline
0.1\%  & 2832  & 1  & 0  & 4 \\ \hline
1\%   &  2942 &  0 & 0 & 4 \\ \hline
10\%   & 2777 &  0 & 0 & 3 \\ \hline
\end{tabular} 		
}
\end{center}
\caption{Comparison of the ILP-3DKP(Z) method for different relative gaps.}
\label{tab:ILPgap}
\end{table}

\noindent The results from Table \ref{tab:ILPgap} indicate that for some instances, a lot of time is spent on improving the solution slightly, or on proving that the optimal solution has been found. This becomes most apparent when increasing the gap from 0.01\% to 0.1\%. In that specific case, the number of nodes decreases from 11750 to \textcolor{black}{1, which illustrates that branching is not required anymore for most instances}. Moreover, from Table \ref{tab:ILPgap}, we can learn that increasing the relative gap from \textcolor{black}{0.1\%} to 10\% has a rather limited effect on the average computation time (i.e., computation time reduces from 4 to 3 seconds, on average). \bigskip

\noindent We believe that for practical applications (e.g., a disaster response mission), it is important to obtain sufficiently good solutions within hours.
Solving our TSS-3DKP with a relative gap of 0.1\% will exactly do that. Consequently, we will apply this relative gap in the remainder of this paper.

\section{To print or not to print?}

In Section \ref{sec:ILP}, we presented a deterministic equivalent of our two-stage stochastic 3D-printing knapsack problem.
We have shown that, with this formulation, it is possible to solve our original problem with standard solvers in reasonable computation time.
Consequently, this formulation can be used to identify whether bringing 3D-printers to a disaster area is useful, or not. As discussed in the introduction of this paper (i.e., in Section 1),
 this decision depends, amongst others, on the following five aspects: $(i)$ the quality of printed items, $(ii)$ the volume and weight of 3D-printers, $(iii)$
 the storage efficiency of printing material, $(iv)$ the \textcolor{black}{3D-}printing time, and $(v)$ demand uncertainty.
In this section, we will study the effect of these five aspects on the outcomes of our TSS-3DKP.
For that, we will, for each aspect, construct a number of instance sets. These instance sets will be constructed such that they vary in the associated aspect only (e.g., a varying quality of the printed items $\alpha$, or a varying weight $w_p$ and volume $v_p$ of the 3D-printer). Per instance set, we then answer the following two questions:

\begin{enumerate}
\item{Should we bring 3D-printers, and if so, how many?}
\item{What is the added value of using 3D-printers?}
\end{enumerate}

\noindent For answering the first question, we will return, per instance set, the median, minimum, maximum, and average number of 3D-printers taken. For answering the second question, we will first determine, per instance, the percentage increase in total reward if bringing 3D-printers is allowed. We then answer our second question by taking the median, minimum, maximum, and average of this percentage over all instances of the instance set. Since we construct several instance sets per aspect, we will answer these questions several times per aspect. We now continue by discussing these outcomes per aspect.

\subsection{The quality of printed items}
\label{sec:alpha}

As discussed in the introduction of this paper, the quality of printed items is still inferior to physical items. This could prevent a decision maker from packing 3D-printers. In this section, we will study the effect of this quality on our two questions. We model the quality of the printed items by solving instance set N100D100S50 for different values of $\alpha$. In particular, we will solve this instance set for $\alpha \in \{0.1,0.2, \ldots,1.0\}$. The results are presented in Table \ref{tab:alpha}. \bigskip

\begin{table}[h!]
\begin{center}
\scalebox{1.0}{
\begin{tabular}{ | l ||  c |  c || c |  c  || }
\hline
 & \multicolumn{2}{ c ||}{\textbf{3D-printers (\#)}} & \multicolumn{2}{ c ||}{\textbf{reward (\%)}}   \\
\hline
$\alpha$ & median & min-max (average)  & median & min-max (average) \\
\hline
0.1	& 0 & 0-0 (0.0) & 0.0 & 0-0 (0.0)	\\ \hline
0.2 & 0  & 0-1 (0.0) & 0.0 & 0-0 (0.0)  \\ \hline
0.3	& 1 & 0-1 (0.5) & 0.1 & 0-3 (0.5)	\\ \hline
0.4	& 1 & 0-1 (0.9) & 2.2 & 0-8 (2.5)	\\ \hline
0.5	& 1 & 1-2 (1.0) & 6.0 & 0-13 (5.7) 	\\ \hline
0.6	& 1 & 1-2 (1.1) & 9.6 & 2-21 (9.7)	\\ \hline
0.7	 & 1 & 1-3 (1.2) & 14.0 & 4-32 (14.3) 	\\ \hline
0.8	 & 1 & 1-3 (1.3) & 19.2 & 7-45 (19.5)	\\ \hline
0.9	 & 1 & 1-3 (1.4)& 25.0 & 9-58 (25.5)	\\ \hline
1.0	&  1 & 1-4 (1.6) & 31.8 & 13-73 (32.5)  \\ \hline
\end{tabular}
}
\end{center}
\caption{An overview of the median, minimum, maximum and average number of both the \textcolor{black}{packed} 3D-printers and the percentage increase in reward, for varying values of $\alpha$.}
\label{tab:alpha}
\end{table}

\noindent From Table \ref{tab:alpha} we can learn that an increase in $\alpha$ (i.e., a higher quality of printed items) coincides with an increase in ($i$) the number of packed 3D-printers and ($ii$) the percentage increase in total reward. Next, we observe that for values of $\alpha$ larger or equal than $30\%$, it is beneficial to bring 3D-printers, implying that printing items, even with a fairly low quality, payoffs rapidly.

\subsection{The effect of the weight and volume of 3D-printers}

\noindent The weight and volume of a 3D-printer is still significant nowadays. This weight and volume can also be used for transporting physical items. In this section, we will study the effect of the 3D-printers weight and volume on our two questions. We will do so by introducing a factor $k$.
This factor represents the number of times a 3D-printer fits (in terms of weight and volume) within the knapsack. In particular, we study instance set N100D100S50 with $k \in \{2,3,5,10,20,30,40,50,\infty\}$.
If $k=\infty$, we set $w_p = v_p = 0$ and otherwise we set $w_p = \nint{\frac{W}{k}}$ and $v_p = \nint{\frac{V}{k}}$.
The results are presented in Table \ref{tab:printer}. \bigskip

\begin{table}[h!]
\begin{center}
\scalebox{1.0}{
\begin{tabular}{ | l ||  c |  c || c |  c  || }
\hline
 & \multicolumn{2}{ c ||}{\textbf{3D-printers (\#)}} & \multicolumn{2}{ c ||}{\textbf{reward (\%)}}   \\
\hline
$k$& median & min-max (average)  & median & min-max (average) \\
\hline
2 &	0 &	0-0 (0.0) &  0 & 0-0 (0.0) \\ \hline
3 & 1 & 0-1	(0.7) & 5.1 & 0-21 (5.7) \\ \hline
5 & 1 & 0-1	(1.0) &	 11.4 &	0-32 (12.2) \\ \hline
10 & 1 & 1-2 (1.1) & 16.1 &	3-39 (16.9) \\ \hline
20 & 1 & 1-2 (1.2)& 18.8 & 4-43	(19.4)\\ \hline
30 & 1 & 1-3 (1.4) &  19.8 & 5-44 (20.3) \\ \hline
40 & 1 & 1-3 (1.4) & 20.3 &	5-45 (20.7) \\ \hline
50 & 1 & 1-3 (1.5) & 20.8 & 5-45 (21.0) \\ \hline
$\infty$ &	3 &	2-7	(3.1) &	22.7 & 6-47	(22.2) \\ \hline
\end{tabular}
}
\end{center}
\caption{\textcolor{black}{An overview of the median, minimum, maximum and average number of both the packed 3D-printers and the percentage increase in reward, for varying values of $k$.}}
\label{tab:printer}
\end{table}

\noindent The results of Table \ref{tab:printer} show that the number of packed 3D-printers and the percentage increase in total reward increases in $k$.
We observe that even for fairly low values of $k$, it is already beneficial to pack a 3D-printer.
The weight and volume of a 3D-printer should exceed 33\% ($k<3$) of the total knapsack weight and volume before packing physical items only (and thus no 3D-printer) is more beneficial.

\subsection{The storage efficiency of printing material}

\noindent As discussed in the introduction of this paper, the storage efficiency of physical items is inferior to the storage efficiency of printing material.
In this section, we study the effect of the storage efficiency of printing material on our two questions.
We do this by introducing factor $l$.
This factor indicates the fraction of printing material needed in comparison to physical items, in terms of both weight as volume.
In particular, we will study instance set N100D100S50 with $l \in \{0,0.1,0.2,...,1\}$ and set $w_i = v_i$, $m_i = \nint{l \cdot w_i}$, and
$w_b = v_b = 1$. The results are presented in Table \ref{tab:mij}. \bigskip

\begin{table}[h!]
\begin{center}
\scalebox{1.0}{
\begin{tabular}{ | l ||  c |  c || c |  c  || }
\hline
 & \multicolumn{2}{ c ||}{\textbf{3D-printers (\#)}} & \multicolumn{2}{ c ||}{\textbf{reward (\%)}}   \\
\hline
$l$ & median & min-max (average)  & median & min-max (average) \\
\hline
$0$ &	1 &	1-3 (1.3) &	17.9 &	6-51 (19.5) \\ \hline
$0.1$ &	1	& 1-3 (1.3) & 16.1 & 6-45 (17.4) \\ \hline
$0.2$ & 1 & 1-2 (1.2) &	14.3 & 5-39	(15.3) \\ \hline
$0.3$ & 1 & 1-2 (1.2) & 12.4 & 5-32	(13.1) \\ \hline
$0.4$ & 1 & 1-2	(1.1) &	10.6 & 4-25	(10.9) \\ \hline
$0.5$ & 1 & 1-2	(1.1) &	8.6 & 3-18 (8.7) \\ \hline
$0.6$ & 1 & 1-1	(1.0) &	6.7	& 3-12 (6.8) \\ \hline
$0.7$ &	1 & 1-1	(1.0) &	5.0 & 2-9 (4.9)\\ \hline
$0.8$ & 1 & 0-1	(1.0) &	3.3 & 0-7 (3.3)\\ \hline
$0.9$ & 1 & 0-1	(0.9) &	2.1 & 0-5 (2.0)\\ \hline
1 &	1 &	0-1	(0.8) &	1.3 & 0-3 (1.2)\\ \hline
\end{tabular}
}
\end{center}
\caption{\textcolor{black}{An overview of the median, minimum, maximum and average number of both the packed 3D-printers and the percentage increase in reward, for varying values of $l$.}}
\label{tab:mij}
\end{table}

\noindent The results of Table \ref{tab:mij} show that the number of 3D-printers and the percentage increase in total reward decreases for an increasing $l$
(i.e. if the material needed to print items takes up more space).
Even when there is no benefit of the storage efficiency of the printing material (i.e., $l=1$) a 3D-printer is still often packed.
This illustrates that the uncertainty \textcolor{black}{in} demand can be enough to bring 3D-printers.

\subsection{3D-\textcolor{black}{p}rinting time}

\noindent Time is crucial during a disaster response mission. However, the time to print a single item may still require hours. This may results in a relative low number of printed items. This number will increase if the time to print a single item reduces (e.g., due to improvements in printing technology). In this section, we will investigate how a change in this printing time effects our two questions. In particular, we will model this by introducing a factor $m \in \{0, 0.1, 0.2,...,1.0, \infty\}$. This factor indicates the fraction of items that can be printed on one 3D-printer, on average. Given this setup, we will use maximum amount of printing time $T = \nint{m \cdot \sum_{s \in S}q_s \sum_{i \in N^p}  t_{i} d_{i}^s}$. We will present our results in Table \ref{tab:T}.\bigskip

\begin{table}[h!]
\begin{center}
\scalebox{1.0}{
\begin{tabular}{ | l ||  c |  c || c |  c  || }
\hline
 & \multicolumn{2}{ c ||}{\textbf{3D-printers (\#)}} & \multicolumn{2}{ c ||}{\textbf{reward (\%)}}   \\
\hline
$m$ & median & min-max (average)  & median & min-max (average) \\
\hline
0 & 0 &	0-0 (0.0) &	0 & 0-0 (0.0) \\ \hline
0.1& 3 & 1-5 (2.8) & 10.3 & 1-21 (10.3) \\ \hline
0.2 & 2& 1-3 (2.2)	& 14.7 & 4-28 (15.1) \\ \hline
0.3 &2& 1-2 (1.8)	& 16.8 & 5-33 (16.9) \\ \hline
0.4 &1	& 1-2 (1.4)	& 17.4 & 5-35 (17.7) \\ \hline
0.5  &1& 1-2 (1.1)	& 18.2 & 5-35 (18.5) \\ \hline
0.6 & 1 & 1-1 (1.0) & 18.7 & 5-37 (19.1) \\ \hline
0.7 & 1 & 1-1 (1.0) & 19.0 & 5-38 (19.4) \\ \hline
0.8	& 1 & 1-1 (1.0) & 19.2 & 5-39 (19.5) \\ \hline
0.9	& 1 & 1-1 (1.0) & 19.2 & 5-39 (19.5) \\ \hline
1 &	1 &	1-1 (1.0) &	19.2 & 5-39 (19.5) \\ \hline
$\infty$ & 1 & 1-1 (1.0) & 19.2 & 5-39 (19.5) \\ \hline
\end{tabular}
}
\end{center}
\caption{\textcolor{black}{An overview of the median, minimum, maximum and average number of both the packed 3D-printers and the percentage increase in reward, for varying values of $m$.}}
\label{tab:T}
\end{table}

\noindent From Table \ref{tab:T}, we observe that, except for the case $m=0$, in which no 3D-printers are taken,
the number of packed 3D-printers decreases in $m$ and converges for large $m$ to exactly one 3D-printer.
Note that for $m=\infty$ all items can be printed on one 3D-printer and consequently there is no reason to pack more than one 3D-printer.
From Table \ref{tab:T}, we also observe that for $0.6 \leq m\leq  1$ at most one 3D-printer is packed,
although not all items can be printed on one 3D-printer. It turns out that in this interval, on average, more than half of the printing time of the demand is utilized on the first 3D-printer. Consequently, the second 3D-printer is, on average, not fully utilized, which makes it likely that this second 3D-printer is not being packed at all. Furthermore, we observe from Table \ref{tab:T} that the percentage increase in total reward increases until $m=0.8$ and stabilizes afterwards. This indicates that for large $m$, some items are preferred being packed physically (e.g. to get a better reward), although these items could be printed without packing any
additional 3D-printers. \bigskip

\subsection{The effect of demand uncertainty}

\noindent One of main advantages of using 3D-printers is the possibility to deal with demand uncertainty. Consequently, it is likely that the number of packed 3D-printers increases for an increasing degree of uncertainty. In this section, we will investigate the effect of this degree of uncertainty on our two questions. We will model the degree of demand uncertainty by solving our instance sets for different values of the upper limit of demand $D$. In particular, we will study instance set N100D100S50 with $D \in \{2^0,2^1,2^2, \ldots, 2^{17}\}$. Recall that, in our instance generator, $D$ identifies the domain of the uniform distributions that generate the demand for all items. Consequently, an increasing $D$ will coincide with a higher degree of uncertainty. Moreover, recall that in our instance generator, the knapsack weight and volume capacity, and the printing time are based on the total average demand. In this section, we are only interested in the effect of changing $D$. For that reason, we fix the following parameters $W=V= 100,000$ and $T=4000$. We report our results in Table \ref{tab:demand}. \bigskip

\begin{table}[h!]
\begin{center}
\scalebox{1.0}{
\begin{tabular}{ | c ||  c |  c || c |  c  || }
\hline
 & \multicolumn{2}{ c ||}{\textbf{3D-printers (\#)}} & \multicolumn{2}{ c ||}{\textbf{reward (\%)}}   \\
\hline
$D$ & median & min-max (average)  & median & min-max (average) \\
\hline
1 & 0 & 0-0 (0.0) & 0.0 & 0-0 (0.0) \\ \hline
2 & 0 & 0-0 (0.0) & 0.0 & 0-0 (0.0) \\ \hline
4 & 0 & 0-0 (0.0) & 0.0 & 0-0 (0.0) \\ \hline
8 & 0 & 0-0 (0.0) & 0.0 & 0-0 (0.0) \\ \hline
16&	1 &	0-1	(0.6) &	0.1 &	0-2 (0.3)\\ \hline
32 & 1 & 1-1 (1.0) & 5.9 & 2-15	(6.4)\\ \hline								
64&	1 &	1-1	(1.0) & 14.7&	6-31 (15.6)\\ \hline
128	& 1	& 1-2 (1.5) & 23.6 & 12-47 (24.3)\\ \hline							
256 & 2 & 1-3	(2.1)&	25.1& 15-51 (26.5)\\ \hline
512	& 3 & 2-4 (2.6) & 22.5 & 8-43 (23.1)\\ \hline		
1024 & 3 & 2-5 (3.1)& 18.0 & 5-60 (20.5)\\ \hline
2048 & 3 & 1-7 (3.3) & 14.5 & 3-37 (15.4)	\\ \hline						
4096 &3 & 0-7 (3.2)	& 10.5 & 0-51 (12.9) \\ \hline
8192 & 3 & 0-8 (3.1) & 6.7 & 0-51 (9.4) \\ \hline			
16384 & 2 & 0-8 (2.4) & 3.0 & 0-67 (7.3)\\ \hline
32768 & 0 & 0-12 (1.5) & 0.0 & 0-55 (4.8) \\ \hline							
65536&	0 & 0-11 (1.0)& 0.0 & 0-31 (2.2)\\ \hline
131072 & 0 &0-11 (1.1) & 0.0 &0-59 (3.5) \\ \hline
							

\end{tabular}
}
\end{center}
\caption{\textcolor{black}{An overview of the median, minimum, maximum and average number of both the packed 3D-printers and the percentage increase in reward, for varying values of $D$.}}
\label{tab:demand}
\end{table}

\noindent From Table \ref{tab:demand}, we observe that both the number of 3D-printers as the percentage increase in reward increase up to a certain degree of uncertainty (i.e., $D=2048$ and $D=256$, respectively), and decrease afterwards. The increasing part confirms our intuition that more uncertainty coincides with packing more 3D-printers. The decreasing parts, however, is counterintuitive. This outcome, which is a result of our way to model uncertainty, can be explained as follows. For sufficiently large degrees of uncertainty, demand of items will exceed the knapsack weight and volume capacity for a majority of the scenarios.
As a consequence, many of these scenarios can be recognized as if there was no restriction on demand for these items.
This means that we can recognize our knapsack problem as a deterministic one. Consequently, there is no reason to pack 3D-printers to adapt for demand uncertainty. The only remaining reason to pack 3D-printers is the better packing efficiency of printing material. Hence, the number of 3D-printers packed will reduce significantly and this is exactly what we see in Table \ref{tab:demand}.

\section{Conclusion}

\noindent In this paper, we shed light on when to bring 3D-printers or not, which amount of printing
material to pack and which items to bring physically, to a disaster area. For that, we introduce a new
type of problem, which we call the two-stage stochastic 3D-printing knapsack problem.  We present a deterministic equivalent of our problem, which allows us to solve the problem with standard solver techniques. Our numerical results illustrates how the quality of printed items, the weight and volume of 3D-printers, the storage efficiency of
printing material, the printing time, and the demand uncertainty effect the outcomes of our model. It turns out that for most situations packing a 3D-printer is beneficial. Only in extreme circumstances, where the quality of printed items is extremely low,
 the size of the 3D-printer is extremely large compared to the knapsack size, when there is no time to print the items, or when demand for items is low, taking no 3D-printers is the best option.\bigskip

\noindent For future research, we identify four interesting research directions. The first one relates to the modelling assumptions of the 3D-printers. In our paper, we consider a sort of "ideal" 3D-printer: it can print all items with the same type of printing material and it can be set-up, and used easily in a disaster response mission. It would be of interest to see how robust our model is against these type of assumptions (e.g., what happens if \textcolor{black}{we have different types of 3D-printers that can print specific types of items only?}).
As a second research direction, one could extend our TSS-3DKP by including dependencies between the items. With this modelling feature, it is possible to include composed items that consist of several printed or physical subitems.
As a third research direction, one could extend our TSS-3DKP by introducing decision dependent uncertainty (cf. \citet{Goel:2006}).
In such a setting, one could invest time to get more accurate demand information. Clearly, this creates a trade-off between reducing demand uncertainty and arriving late(r) at a disaster area. As a final research direction, it is of interest to see how our TSS-3DKP could be applied to other application domains (e.g., a military or aerospace setting), which each may have their own, unique, characteristics. \bigskip

\noindent One of the strengths of our model is that sufficiently large instances can be solved by commercial solvers within reasonable computation time.
However, this may no longer be the case for (the aforementioned) future research directions. Consequently, the need for (other) efficient solution methods becomes apparent.
A candidate for this is logic-based Benders decomposition \citep{Hooker:2003}, a generalization of Benders
decomposition \citep{Benders:1962} that allows the second-stage problem to be any optimization problem rather than precisely a linear or nonlinear programming problem.

\bibliography{reference}
\bibliographystyle{abbrvnat}

\end{document}